\documentclass[11pt,a4paper,leqno]{article}
\pdfoutput=1

\usepackage[latin1]{inputenc}
\usepackage[english]{babel}
\usepackage[margin=2.75cm]{geometry}

\usepackage{authblk}

\usepackage[sort&compress,numbers]{natbib}
\usepackage[titletoc,title]{appendix}

\usepackage{amsmath}
\usepackage{amsfonts}
\usepackage{amssymb}
\usepackage{amsthm}
\usepackage{bbm}									
\usepackage{mathtools}								
\usepackage{xcolor}
\usepackage{tikz}

\usepackage{graphicx}
\usepackage[suffix=]{epstopdf}
\usepackage{booktabs}
\usepackage[font=small,labelfont=bf]{caption} 

\definecolor{thelinkcolor}{RGB}{0,0,150}
\usepackage{hyperref}
\usepackage[nameinlink]{cleveref}
\crefname{subsection}{subsection}{subsections}

\crefname{equation}{}{}
\crefname{theorem}{theorem}{theorems}
\crefname{example}{example}{examples}
\crefname{lemma}{lemma}{lemmas}
\crefname{proposition}{proposition}{propositions}
\crefname{figure}{figure}{figures}
\crefname{table}{table}{tables}
\crefname{subsection}{subsection}{subsections}

\Crefname{equation}{}{}
\Crefname{theorem}{Theorem}{Theorems}
\Crefname{example}{Example}{Examples}
\Crefname{lemma}{Lemma}{Lemma}
\Crefname{proposition}{Proposition}{Proposition}
\Crefname{figure}{Figure}{Figures}
\Crefname{table}{Table}{Tables}
\Crefname{subsection}{Subsection}{Subsections}

\hypersetup{colorlinks=true,
	urlcolor=blue,
	linkcolor=blue,
	citecolor=blue,
	bookmarksdepth=paragraph,
	}

\usepackage[shortlabels]{enumitem}
\setlist[itemize]{itemindent=0ex,itemsep=-0.5ex,leftmargin=2ex,topsep=5pt}
\setlist[enumerate]{label={\arabic*)},leftmargin=*,nolistsep,noitemsep,topsep=1ex}

\newcommand{\R}{\mathbb{R}}							

\renewcommand{\vec}[1]{\mathbf{#1}}    	  			
\newcommand{\Rey}{\mbox{\it Re}} 		  			
\newcommand{\dt}{\,{\rm d}t} 			  		    
\newcommand{\ssep}{\,;}					  		    
\newcommand{\abs}[1]{\left\vert #1 \right\vert}		

\makeatletter
\newcommand{\subalign}[1]{%
	\vcenter{%
		\Let@ \restore@math@cr \default@tag
		\baselineskip\fontdimen10 \scriptfont\tw@
		\advance\baselineskip\fontdimen12 \scriptfont\tw@
		\lineskip\thr@@\fontdimen8 \scriptfont\thr@@
		\lineskiplimit\lineskip
		\ialign{\hfil$\m@th\scriptstyle##$&$\m@th\scriptstyle{}##$\crcr
			#1\crcr
		}%
	}
}
\makeatother

\definecolor{RGBblue}{rgb}{0,0,1}
\definecolor{RGBred}{rgb}{1,0,0}
\definecolor{RGBgreen}{rgb}{0,1,0}
\definecolor{RGBmagenta}{rgb}{1,0,1}
\definecolor{RGBblack}{rgb}{0,0,0}

\definecolor{matlabblue}{rgb}{0,0.4470,0.7410}
\definecolor{matlabred}{rgb}{0.6350,0.0780,0.1840}
\definecolor{matlabgreen}{rgb}{0.4660,0.6740,0.1880}
\definecolor{matlaborange}{rgb}{0.8500,0.3250,0.0980}
\definecolor{matlabmagenta}{rgb}{0.4940,0.1840,0.5560}
\definecolor{matlabgrey}{rgb}{0.65,0.65,0.65}

\newcommand\solidrule[1][10pt]{\rule[0.5ex]{#1}{1.5pt}}

\newcommand\dotdashedrule{\mbox{%
		\solidrule[3pt]\hspace{1.5pt}\solidrule[1pt]\hspace{1.5pt}\solidrule[3pt]}}

\newcommand{\mydot}[1]{%
	\protect\begin{tikzpicture}%
	\protect\filldraw[color=#1] (1.5ex,1.5ex) circle (1.8pt);
	\protect\end{tikzpicture}%
}
\newcommand{\mysolidsquare}[1]{%
	\protect\begin{tikzpicture}%
	\protect\filldraw[thick,color=#1] (0,0) -- (0.75ex,0) -- (0.75ex,0.75ex) -- (0,0.75ex) -- (0,0);
	\protect\end{tikzpicture}%
}

\newtheorem{proposition}{Proposition}

\theoremstyle{definition}

\author[1]{Mayur V. Lakshmi\thanks{Email address for correspondence: \href{mailto:mayur.lakshmi11@imperial.ac.uk}{mayur.lakshmi11@imperial.ac.uk}}} 
\author[1]{Giovanni Fantuzzi}
\author[1]{Sergei I. Chernyshenko} 
\author[2]{\\Davide Lasagna}
\affil[1]{Department of Aeronautics, Imperial College London, London, SW7 2AZ, UK}
\affil[2]{Faculty of Engineering and the Environment, University of Southampton, Highfield, Southampton, SO17 1BJ, UK}
\title{\textbf{Finding unstable periodic orbits: a hybrid approach with polynomial optimization}}

\begin{document}
\maketitle
\begin{center}
\begin{minipage}{0.85\linewidth}
\begin{small}
\noindent
\textbf{Abstract.}\enspace\
We present a novel method to compute unstable periodic orbits (UPOs) that optimize the infinite-time average of a given quantity for polynomial ODE systems. The UPO search procedure relies on polynomial optimization to construct nonnegative polynomials whose sublevel sets approximately localize  parts of the optimal UPO, and that can be used to implement a simple yet effective control strategy to reduce the UPO's instability. Precisely, we construct a family of controlled ODE systems parameterized by a scalar $k$ such that the original ODE system is recovered for $k = 0$, and such that the optimal orbit is less unstable, or even stabilized, for $k>0$. Periodic orbits for the controlled system can be more easily converged with traditional methods and numerical continuation in $k$ allows one to recover optimal UPOs for the original system. The effectiveness of this approach is illustrated on three low-dimensional ODE systems with chaotic dynamics.
\end{small}
\end{minipage}
\end{center}
\vskip4ex

\section{Introduction}
\label{s:introduction}
Computing unstable periodic orbits (UPOs) for systems governed by ordinary differential equations (ODEs) is a fundamental problem in the study of nonlinear dynamical systems that exhibit chaotic dynamics. UPOs embedded in a chaotic attractor provide a ``skeleton'' around which trajectories evolve through a continuous process of attraction and repulsion along the orbits' stable and unstable manifolds~\cite{Cvitanovic1991}. UPOs are also fundamental to the periodic orbit theory introduced by Cvitanovi\'{c} \textit{et al.}~\cite{Artuso1990,Cvitanovic1995}, which states that the infinite-time average of any observable of interest over a chaotic trajectory can be expressed as a weighted sum of the time averages over individual UPOs.

A widespread and very effective strategy to find dynamically relevant UPOs is to perform recurrence analysis on long chaotic trajectories in order to identify nearly periodic segments that can be used as initial conditions for a variety of UPO-finding algorithms. The simplest family of such algorithms are shooting methods, which apply the Newton--Raphson algorithm to the Poincar\'{e} return map of the dynamical system to compute one (single shooting) or more (multiple shooting) points on the periodic orbit as well as its period (see, for instance,~\cite{Kelley2003,Ascher1995a} and~\cite[Chapter 12]{Cvitanovic2017}). A more robust family of methods are variational ones~\cite{Lan2004,Boghosian2011a}, which attempt to deform a (discrete) closed loop in state space into a UPO by minimizing a nonconvex cost function that, loosely speaking, measures the deviation between the tangent to the loop and the direction of the ODE's vector field at each point.

Even though this approach to finding UPOs has been applied very successfully to a wide variety of high-dimensional systems, including turbulent fluid flows (see, e.g.,~\cite{Fazendeiro2010,Chandler2013,Lucas2014}), it suffers from an inherent drawback: success relies on the availability of a very good initial approximation to a periodic orbit. In the case of shooting methods this is because trajectories starting from inaccurate initial conditions diverge very quickly from the UPO, whereas the basin of attraction for the Newton--Raphson root finder is often very small~\cite{Kelley2003}. For variational methods, instead, poor initial loops may converge to local minima of the nonconvex cost function that are not periodic orbits of the underlying ODE. Since extracting good initial approximations with recurrence analysis relies on chaotic trajectories shadowing a UPO sufficiently closely for one entire period, one generally can find only orbits that are short and/or have only a few unstable directions. UPOs that are very unstable or lie in parts of the state space rarely visited by chaotic trajectories, however, often remain undetected.

In this work, we present a new strategy for constructing good initial approximations to UPOs that may otherwise be difficult to find. This strategy can be implemented computationally on low-dimensional ODE systems with polynomial dynamics and is part of a broader framework to characterize trajectories that maximize or minimize the infinite-time average of a given quantity of interest $\Phi$. Such \textit{extremal} orbits are useful for control purposes, as their knowledge facilitates the design of optimal control actions to stabilize desired dynamics or suppress undesired ones. Of course, focusing on UPOs that optimize time averages is not restrictive because every periodic orbit is extremal for at least one choice of $\Phi$, namely, the indicator function of the orbit itself. Although this is unknown in practice, varying $\Phi$ potentially enables one to identify a large number of periodic orbits.

Underpinning our new approach is the recent realization that nearly sharp bounds on extreme values of time averages can be derived by constructing so-called auxiliary functions of the system's state~\cite{Chernyshenko2014,Fantuzzi2016,Goluskin2018,Tobasco2018,Korda2021}, which are similar to the Lyapunov functions used in stability analysis but satisfy a different set of constraints. In addition to providing bounds on extremal time averages, Tobasco \textit{et al.}~\cite{Tobasco2018} showed that auxiliary functions provide information about the location and shape of the corresponding extremal trajectory---which, very often, is a UPO~\cite{Yang2000}. These observations have already been exploited by some of the authors~\cite{Lakshmi2020} to compute near-extremal UPOs for a nine-dimensional ODE model of shear flow with single shooting methods. However, the inherent ill-conditioning of initial value problems for chaotic ODEs limits the applicability of this approach to cases in which very nearly optimal auxiliary functions can be constructed accurately, which is often not the case. 

In this work, we go one step further and show that auxiliary functions can also be used to construct an effective open-loop control action that can be expected to reduce the instability of the extremal UPO and, in some cases, provably stabilizes it. Adding this control to the ODE system produces a family of controlled systems parameterized by the control amplitude, to which traditional shooting or variational methods are more easily applied. One can then attempt to numerically continue any orbit computed with control by decreasing the control amplitude, until a UPO for the original ODE is obtained. In particular, combining this approach with the techniques developed in~\cite{Lakshmi2020} enables one to search for extremal UPOs robustly with auxiliary functions that are not sufficiently close to being optimal for the latter to work in isolation.

The rest of the paper is organized as follows. \Cref{s:auxiliary-function-method} reviews the auxiliary function method and how auxiliary functions can be leveraged to localize extremal trajectories. \Cref{s:control-methodology} introduces our control strategy to stabilize UPOs and describes how the construction of a family of controlled systems allows for the computation of extremal UPOs for polynomial ODEs. The practical potential of this approach is demonstrated in \cref{s:results} on three low-dimensional ODE systems that display chaotic dynamics. Computational merits and limitations of our method are discussed in \cref{s:discussion} along with possible directions for further improvement. Finally, \cref{s:conclusion} offers concluding remarks.

\section{Approximating extremal UPOs with auxiliary functions}
\label{s:auxiliary-function-method}
Consider an autonomous dynamical system governed by the ODE
\begin{equation}
\label{e:ode}
\frac{\mathrm{d} \vec{a}}{\dt} = \vec{f}(\vec{a}), \quad \vec{a}(0) = \vec{a}_0,
\end{equation}
where $\vec{a} \in \R^n$ and $\vec{f}: \R^n \to \R^n$ is smooth. The infinite-time average of a function $\Phi(\vec{a})$ along the trajectory starting from $\vec{a}_0$ is defined as
\begin{equation}
\label{e:phi-time-average-definition}
\overline{\Phi} \left(\vec{a}_0\right) := \lim_{\tau \to \infty} \frac{1}{\tau} \int_0^\tau \Phi \left[ \vec{a}\left(t \ssep \vec{a}_0\right) \right] \dt,
\end{equation}
where $\vec{a}(t \ssep \vec{a}_0)$ denotes the trajectory $\vec{a}(t)$ with initial condition $\vec{a}_0$, and we assume for simplicity that the limit exists. We are interested in the maximal value of $\overline{\Phi}$ over all bounded trajectories,
\begin{equation}
\label{e:max-phi-definition}
\overline{\Phi}^* := \max _{\substack{\vec{a}_0 \in \mathbb{R}^n:\\\|\vec{a}(t \ssep \vec{a}_0)\|<\infty \,\forall t}} \overline{\Phi}(\vec{a}_0),
\end{equation}
as well as the initial conditions and corresponding trajectories which achieve it. Note that considering maximal time averages only is not restrictive, as minimizers of $\overline{\Phi}$ coincide with maximizers of $\overline{-\Phi}$.

Upper bounds on $\overline{\Phi}^*$ can be computed in a relatively straightforward way. Suppose that there exist a function $P(\vec{a})$ such that $\overline{P} = 0$ and a constant $U$ such that $\Phi(\vec{a}) + P(\vec{a}) \leq U$ for any $\vec{a}$.  Averaging this inequality along the trajectory $\vec{a}(t \ssep \vec{a}_0)$ yields $\overline{\Phi}(\vec{a}_0) \leq U$ for any $\vec{a}_0$, so $\overline{\Phi}^*\leq U$. To construct a function $P(\vec{a})$ with zero average, one can take $P(\vec{a}) = \vec{f}(\vec{a}) \cdot \nabla V(\vec{a})$ with any $V:\R^n \to \R$ in the class $C^1$ of continuously differentiable functions. Indeed, along any bounded trajectory of \cref{e:ode} the chain rule gives
\begin{align}
\overline{P[\vec{a}(t \ssep \vec{a}_0)]}
&= \overline{\vec{f}[\vec{a}(t \ssep \vec{a}_0)] \cdot \nabla V[\vec{a}(t \ssep \vec{a}_0)]} \\[1ex]
&= \overline{\frac{\mathrm{d}}{\dt} V[\vec{a}(t \ssep \vec{a}_0)]} \nonumber\\
&= \lim_{\tau \to \infty} \frac{V[\vec{a}(\tau \ssep \vec{a}_0)] - V(\vec{a}_0)}{\tau} \nonumber\\[1ex]
&= 0. \nonumber
\end{align}
The best upper bound on $\overline{\Phi}^*$ is obtained by optimizing over the choice of $V$, hereafter called \textit{auxiliary function}:
\begin{equation}
\label{e:phi-bar-star-bound}
\overline{\Phi}^* \leq \inf_{V \in C^1} \left\{U \mid   U - \Phi(\vec{a}) - \vec{f}(\vec{a}) \cdot \nabla V(\vec{a}) \geq 0 \quad \forall \vec{a} \in \mathbb{R}^n \right\}.
\end{equation}

If $\Phi$ and $\vec{f}$ are polynomial, feasible polynomial auxiliary functions and their corresponding bounds on $\overline{\Phi}^*$ can be constructed computationally upon replacing the inequality constraint in~\cref{e:phi-bar-star-bound} with the stronger requirement that the polynomial $U - \Phi(\vec{a}) - \vec{f}(\vec{a}) \cdot \nabla V(\vec{a})$ be a sum of squares (SOS)~\cite{Chernyshenko2014,Fantuzzi2016,Goluskin2018,Goluskin2019,Lasagna2016,Lakshmi2020}, which can be handled using efficient algorithms for convex optimization~\cite{Parrilo2003,Parrilo2012,Lasserre2015}. Moreover, if all trajectories of~\cref{e:ode} are absorbed by a compact set $\Omega$, then the constraint in~\cref{e:phi-bar-star-bound} can be restricted to $\Omega$ and the corresponding SOS computations are guaranteed to return arbitrarily sharp bounds on $\overline{\Phi}^*$~\cite{Tobasco2018,Lakshmi2020}. More precisely, for any $\delta>0$ one can construct a polynomial auxiliary function $V_{\delta}$ which provides a bound $U_{\delta}$ with $\overline{\Phi}^* \leq U _{\delta} \leq \overline{\Phi}^* + \ \delta$. Importantly, this often happens in practice even when no absorbing set $\Omega$ exists.

The crucial observation for this work is that any near-optimal auxiliary function constructed with polynomial optimization not only produces an upper bound on $\overline{\Phi}^*$, but can also be used to localize the associated extremal trajectories in state space -- which are, very often, UPOs. To understand this localization framework, suppose first that there exists an optimal auxiliary function $V_0$ proving a sharp upper bound $U_0=\overline{\Phi}^*$. Then, any extremal trajectory $\vec{a}(t)$ must satisfy~\cite{Fantuzzi2016,Tobasco2018}
\begin{equation}
\label{e:inf-time-average-optimal-V}
\overline{\overline{\Phi}^* - \Phi[\vec{a}(t)] -  \vec{f}[\vec{a}(t)] \cdot \nabla V_0[\vec{a}(t)]}= 0.
\end{equation}
Since the quantity being averaged in~\cref{e:inf-time-average-optimal-V} is nonnegative it must actually vanish pointwise in time, so the extremal trajectory necessarily lies inside the set
\begin{equation}
\mathcal{S}_0 := \{ \vec{a} \in \mathbb{R}^n \mid \overline{\Phi}^* - \left[ \vec{f}(\vec{a}) \cdot \nabla V_0(\vec{a}) + \Phi(\vec{a}) \right] = 0 \}.
\end{equation}
Although this set may contain also points that are not on the extremal periodic orbit, it provides guidance to locate the extremal trajectory. 

In practice, an optimal auxiliary function $V_0$ is rarely available and one can only construct a $\delta$-suboptimal one, $V_\delta$, with corresponding bound $U_\delta$. Nevertheless, it was shown in~\cite{Tobasco2018} that, if the extremal trajectory is a periodic orbit, then for any $\varepsilon \geq \delta$ it must lie inside the set
\begin{equation}
\label{e:tobasco-set}
\mathcal{S}_{\varepsilon} = \{ \vec{a} \in \mathbb{R}^n \mid D(\vec{a}) \leq \varepsilon \}
\end{equation}
for a fraction of its time period no smaller than $F := 1 - \delta/\varepsilon$, where
\begin{equation}
\label{e:polynomial}
D(\vec{a}) := U_{\delta} -  \vec{f}(\vec{a}) \cdot \nabla V_{\delta}(\vec{a}) -\Phi(\vec{a}).
\end{equation}
Since the polynomial $D(\vec{a})$ is nonnegative by construction, when $\delta \ll 1$ (that is, when the auxiliary function is close to being optimal) taking $\delta \ll \varepsilon \ll 1$ often allows one to keep $F\approx 1$ while excluding much of the state space from $\mathcal{S}_{\varepsilon}$, which therefore localizes large portions of extremal and near-extremal orbits.

\begin{figure}
\centering
\includegraphics[width=0.9\textwidth]{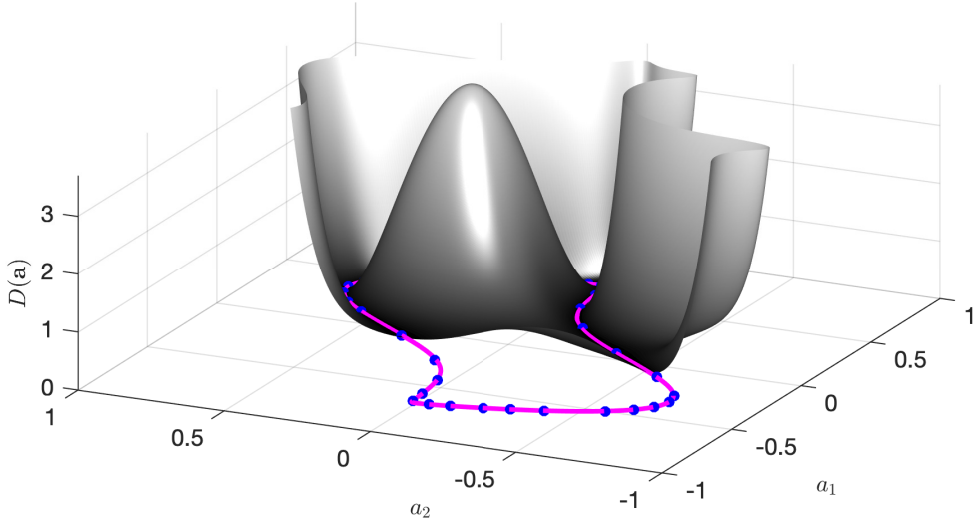}
\caption{Surface plot of the polynomial $D(\vec{a})$ for the rescaled van der Pol oscillator~\cref{e:VdP-system} corresponding to the observable $\Phi(\vec{a}) = a_1^2 + a_2^2$ and a polynomial auxiliary function $V$ of degree 16. The solid line ({\color{RGBmagenta}\solidrule})  indicates the limit cycle for~\cref{e:VdP-system}, which is the extremal trajectory for the chosen $\Phi$. Blue dots (\mydot{RGBblue}) show points lying in $\mathcal{S}_{10^{-4}}$ obtained via direct unconstrained minimization of $D(\vec{a})$.}
\label{fig:VdP-surface-plot}
\end{figure}

While constructing the entire set $\mathcal{S}_{\varepsilon}$ is computationally intractable except for ODE systems of very low dimension, obtaining points that lie in it is relatively straightforward. For instance, global minimizers for $D(\vec{a})$ lie in $\mathcal{S}_{\varepsilon}$ for all $\varepsilon \geq \delta$ and can sometimes be recovered directly from the solution of the optimization problem resulting from an SOS relaxation of~\cref{e:phi-bar-star-bound}~\cite[section 6.1.2]{Lasserre2015}. A simpler and more robust procedure, however, is to numerically search for points where $D(\vec{a}) \leq \varepsilon$  by minimizing $D$ using any nonlinear minimization algorithm, initialized from a large number of initial conditions~\cite{Lakshmi2020}. Since all points in $\mathcal{S}_{\varepsilon}$ along the extremal periodic orbit are close to being global minimizers for $D$ when $\varepsilon$ is small, it is reasonable to expect that $\nabla D$ will be small along the part of the extremal periodic orbit contained in $\mathcal{S}_{\varepsilon}$ and large elsewhere. The minimization routine will therefore quickly descend to this flat region and then slowly approach a local minimum, producing a collection of points in $\mathcal{S}_\varepsilon$ on or close to the extremal periodic orbit as part of the process. This typical situation is illustrated in~\Cref{fig:VdP-surface-plot}, which shows the unstable limit cycle of the reverse-time van der Pol oscillator
\begin{equation}
\label{e:VdP-system}
\frac{\mathrm{d}a_1}{\dt} =-3a_2, \qquad \frac{\mathrm{d}a_2}{\dt} =-4(1-9a_1^2)a_2
\end{equation}
along with the polynomial $D$ obtained with $\Phi(\vec{a}) = a_1^2 + a_2^2$ (whose infinite-time average on bounded trajectories is maximized on the limit cycle) and a near-optimal polynomial auxiliary function of degree 16. All points on the limit cycle satisfy $D(\vec{a}) \leq 10^{-6}$.

As demonstrated in~\cite{Lakshmi2020}, any of the points computed by minimizing $D(\vec{a})$ can be used as initial conditions for algorithms that converge to UPOs by evolving the system's dynamics forward in time. However, this basic strategy suffers from two fundamental limitations. The first is that only a finite number of local minima of $D$ may exist when $V$ is suboptimal, to which the minimization routine used will repeatedly converge. The resulting approximation of the extremal periodic orbit can therefore be sparse, making the use of a multiple-shooting or variational method to converge the UPO infeasible in practice. The second fundamental limitation is that, for a given near-optimal auxiliary function and a given $\varepsilon$, the set $\mathcal{S}_{\varepsilon}$ contains points not on the extremal UPO and possibly far from it~\cite[section 4]{Tobasco2018}. Thus, the point obtained via minimization of $D(\vec{a})$ may not be a good initial condition for single-shooting techniques. In the next section we address the first of these limitations by introducing a novel control methodology to reduce the instability of extremal UPOs.

\section{Control methodology}
\label{s:control-methodology}
Auxiliary functions can be used not only to localize extremal trajectories, but also be used to formulate an effective control strategy to stabilize them or, at least, reduce their instability. This observation, which is the main contribution of this work, can be used to aid the computation of extremal UPOs with traditional shooting or variational techniques. \Cref{s:control-optimal} describes this control strategy in the case of optimal auxiliary functions, for which a rigorous stabilization result can be established without much difficulty. The suboptimal case, which is more relevant in practice but for which we do not have similar theoretical results, is discussed in \cref{s:control-suboptimal}. We shall assume throughout that, as is often the case for ODE systems with chaotic dynamics, the trajectory achieving the maximum time average $\overline{\Phi}^*$ is a UPO $\mathcal{O} := \{\vec{a}(t) \mid 0 \leq t < T\}$ with period $T$.

\subsection{Optimal auxiliary functions}
\label{s:control-optimal}
Suppose that $V$ is an optimal auxiliary function and recall from \cref{s:auxiliary-function-method} that the corresponding nonnegative polynomial $D(\vec{a})$ defined in~\cref{e:polynomial} must vanish everywhere along the extremal UPO $\mathcal{O}$, so $\nabla D = 0$ on it. Furthermore, it is not unreasonable to expect that $D(\vec{a})$ is steep elsewhere, so $\nabla D \neq \vec{0}$ in a neighborhood of $\mathcal{O}$ (cf. \Cref{fig:VdP-surface-plot}). Then, consider the ODE
\begin{equation}
\label{e:controlled-ode}
\frac{\mathrm{d}\vec{a}}{\dt} = \vec{f}(\vec{a}) -k \nabla D(\vec{a}), \quad \vec{a}(0) =\vec{a}_0,
\end{equation}
where $k$ is an arbitrary nonnegative scalar which takes on the role of a control parameter. Since $\nabla D =0 $ on the extremal UPO $\mathcal{O}$ by construction, $\mathcal{O}$ remains a periodic orbit of~\cref{e:controlled-ode} for any choice of $k$, but its stability properties depend crucially on $k$.

To see this, observe that the two terms on the right-hand side of~\cref{e:controlled-ode} have a clear dynamical meaning. Given an initial condition $\vec{a}_0$ sufficiently close to (but not on) the extremal UPO $\mathcal{O}$, the vector field $\vec{f}$ approximately advances the ensuing trajectory along the periodic orbit, but its instability eventually leads to divergence. The $-k\nabla D$ term counteracts this instability by pushing the trajectory back towards the set $\mathcal{S}_0 = \{\vec{a} \mid D(\vec{a}) = 0\}$, which contains the entire extremal UPO because the auxiliary function $V$ used to construct $D$ is optimal. Thus, the term $-k \nabla D$ effectively acts as a control term that reduces the instability of the UPO and whose authority is proportional to $k$. In particular, the controlled ODE~\cref{e:controlled-ode} reduces to the original ODE~\cref{e:ode} when $k = 0$, while letting $k \to \infty$ and rescaling time by $k^{-1}$ leads to the ODE governing the steepest-descent minimization of~$D$.

With this intuition in mind, it is not difficult to prove that the orbit $\mathcal{O}$ is locally stable for all sufficiently large $k$ provided that $\nabla D$ does not vanish in a neighborhood of $\mathcal{O}$.
\begin{proposition}
\label{p:proposition}
Assume that the ODE~\cref{e:ode} has a periodic orbit $\mathcal{O} = \{ \vec{a}(t) \}_{t \in (0,T]}$. Suppose also that $D(\vec{a})$ is continuously differentiable, that $D(\vec{a}) = 0$ on $\mathcal{O}$, and that $D(\vec{a}) \geq 0$ in a set $\Omega \subseteq \R^n$ that contains $\mathcal{O}$. Further, assume that there exists a bounded open neighborhood $\mathcal{N}$ of the orbit $\mathcal{O}$ such that $\| \nabla D(\vec{a}) \| > 0$ for all $\vec{a} \in \mathcal{N} \setminus\mathcal{O}$. Then, there exists $k_0$ such that $\mathcal{O}$ is a locally stable orbit for the controlled system~\cref{e:controlled-ode} for all $k > k_0$.
\begin{proof}
Since $D$ is continuous, there exists $\gamma_0 > 0$ such that the set $\mathcal{U}_{\gamma} = \{ \vec{a} \in \mathcal{N} \mid D(\vec{a}) \leq \gamma\}$ is compact for all $\gamma \leq \gamma_0$, so it does not intersect the boundary of $\mathcal{N}$. The same is true for $\mathcal{V}_{\gamma} = \{ \vec{a} \in \mathcal{U}_{\gamma} \mid \frac{1}{2} \gamma \leq D(\vec{a}) \leq \gamma\}$ and
\begin{equation}
c := \min_{\vec{a} \in \mathcal{V}_{\gamma}} \| \nabla D \|^2 > 0
\end{equation}
because $\| \nabla D(\vec{a}) \| > 0$ on $\mathcal{V}_{\gamma}$.
Then,
\begin{equation}
\frac{\mathrm{d}}{\dt} D(\vec{a}(t)) 
= \vec{f} \cdot \nabla D - k \| \nabla D \|^2 
\leq \max_{\vec{a}\in\mathcal{V}_\gamma}\| \vec{f} \cdot \nabla D\| - kc
\end{equation}
along any trajectory of~\cref{e:controlled-ode} starting in $\mathcal{V}_{\gamma}$.
Consequently, $D$ decays along trajectories provided that
\begin{equation}
k > \max_{\vec{a}\in\mathcal{V}_\gamma}\frac{\| \vec{f} \cdot \nabla D\|}{c}  =: k_0.
\end{equation}
In particular, for all such $k$ trajectories cannot escape the set $\mathcal{U}_{\gamma}$, which is a neighborhood of the orbit $\mathcal{O}$. Moreover, since $\mathcal{O} = \bigcap_{\gamma > 0} \mathcal{U}_{\gamma}$ we can make this trapping set $\mathcal{U}_{\gamma}$ arbitrarily small by taking $\gamma$ arbitrarily small. This proves that $\mathcal{O}$ is locally stable, as claimed.
\end{proof}
\end{proposition}

In principle, this stability result enables the computation of the UPO $\mathcal{O}$ for system~\eqref{e:ode} in a straightforward way: simply find points $\vec{a}_0$ with $D(\vec{a}_0) \leq \varepsilon$ using nonlinear minimization algorithms, as described in \cref{s:auxiliary-function-method}, and use them as initial conditions to simulate the controlled ODE~\cref{e:controlled-ode}. The process can be repeated with increasingly small $\varepsilon$ and increasingly large $k$ until $\vec{a}_0$ falls within the basin of attraction of $\mathcal{O}$. In practice, however, optimal auxiliary functions required by this approach are rarely available in practice, so one must adjust the procedure to allow for suboptimal ones. We turn to this next.

\subsection{Suboptimal auxiliary functions}
\label{s:control-suboptimal}
The controlled ODE~\cref{e:controlled-ode} can be formulated and solved numerically even when the auxiliary function $V$ used to construct $D$ is suboptimal, but two complications arise. First, the extremal orbit $\mathcal{O}$ is generally not a trajectory of the controlled system, because $D$ and $\nabla D$ need not vanish along it when $V$ is suboptimal. Second, and most important, it is possible that taking $k\gg1$ introduces unwanted stable equilibria in the vicinity of $\mathcal{O}$ where $\vec{f}-k\nabla D$ vanishes, preventing trajectories of the controlled system from shadowing $\mathcal{O}$ over an entire period.

To avoid this issue and increase the likelihood that the controlled system possesses a periodic orbit that continuously deforms into $\mathcal{O}$ as $k$ is reduced, we project the control term $-k \nabla D$ onto the subspace perpendicular to $\vec{f}$ and replace~\cref{e:controlled-ode} with
\begin{equation}
\label{e:controlled-projection-ode}
\frac{\mathrm{d}\vec{a}}{\dt} = \vec{f}(\vec{a}) -k \left[ I - \frac{\vec{f}\otimes\vec{f}}{\|\vec{f}\|^2} \right] \nabla D(\vec{a}), 
=: \vec{h}_k(\vec{a}), 
\qquad \vec{a}(0) = \vec{a}_0.
\end{equation}

As proved in~\Cref{app:stability-result}, the stability result in~\Cref{p:proposition} for optimal $V$ extends to this modified controlled system under moderate assumptions on the behavior of $D$ and $\vec{f}$ near the extremal orbit $\mathcal{O}$. For suboptimal $V$, however, one cannot guarantee that increasing $k$ will result in the existence of a stable periodic orbit for~\cref{e:controlled-projection-ode}. Nevertheless, provided that $V$ is sufficiently close to optimal, it is not unreasonable to expect that there exists a family of UPOs for the family of systems~\cref{e:controlled-projection-ode} that connects to the extremal UPO for the original system~\cref{e:ode}. In addition, if $D$ increases rapidly in the directions normal to $\vec{f}$ in the vicinity of the extremal UPO $\mathcal{O}$ as in~\Cref{fig:VdP-surface-plot}, periodic orbits for large $k$ are likely to be less unstable than the original extremal orbit $\mathcal{O}$ obtained with $k=0$ because the control term still strongly damps perturbations normal to it.

These heuristic observations suggest that traditional techniques for computing UPOs may be much more effective for large values of $k$ than for $k=0$. Moreover, since at least part of the original orbit $\mathcal{O}$ must lie in the set $\mathcal{S}_\varepsilon$ defined in~\cref{e:tobasco-set} when $V$ is near-optimal, points in it are not unlikely to be good initial conditions to find exactly the branch of UPOs that connects to $\mathcal{O}$ as $k \to 0$ (if one exists). Therefore, we propose to search for extremal UPOs using the following 4-step procedure:
\begin{enumerate}[(1)]
    \item Construct a near-optimal polynomial auxiliary function $V$ by solving an SOS relaxation of~\cref{e:phi-bar-star-bound} as described in~\cite{Chernyshenko2014,Fantuzzi2016,Goluskin2018,Goluskin2019,Lakshmi2020}.
    \item Construct the polynomial $D$ and, given a tolerance $\varepsilon>0$, attempt to find points where $D(\vec{a})\leq \varepsilon$ with direct nonlinear minimization of $D$.
    \item Fix $k>0$ and search for UPOs of the controlled ODE system~\cref{e:controlled-projection-ode} using the points obtained at step 2 as initial conditions for shooting methods or for recurrence analysis.
    \item Repeat steps 2 and 3 with increasingly small $\varepsilon$ and $k$ until a periodic orbit is found, then perform continuation in $k$ until $k=0$.
\end{enumerate}
We stress that this procedure is not guaranteed to work because the polynomial $D$ may not behave as illustrated in \Cref{fig:VdP-surface-plot}, UPOs for~\cref{e:controlled-projection-ode} may be hard to find even with $k \gg 1$, and any branch of UPOs one manages to find may not continue up to $k=0$. Nevertheless, numerical experiments reveal that our strategy is often successful in practice.

\section{Examples}
\label{s:results}
We now demonstrate the potential of the control methodology described in \cref{s:control-suboptimal} to find extremal UPOs on three low-dimensional ODE systems that display chaotic dynamics. Details of our numerical implementation are discussed in~\Cref{ss:numerical-implementation}.

\subsection{A three-dimensional chaotic system}
\label{ss:sprF-system}
Consider the three-dimensional polynomial ODE system
\begin{equation}
\label{e:sprF-system}
\frac{\mathrm{d} a_1}{\dt} = a_2 + a_3, \qquad \frac{\mathrm{d} a_2}{\dt} = -a_1 + \frac12 a_2, \qquad \frac{\mathrm{d} a_3}{\dt} = a_1^2 - a_3,
\end{equation}
which has two equilibrium points at $(0,0,0)$ and $(-2,-4,4)$  and a chaotic attractor~\cite{Sprott1994}. Trajectories starting outside the basin of the attraction of these invariant structures may become unbounded, but we can still search for extremal UPOs as long as near-optimal auxiliary functions to bound time averages on bounded trajectories can be constructed. This is not guaranteed by the theoretical results in~\cite{Tobasco2018,Lakshmi2020} because~\cref{e:sprF-system} has no compact absorbing set, but it appears to be true in practice.

We therefore applied our 4-step control strategy to search for extremal UPOs that maximize the infinite-time average of the following observables:
\begin{equation}
\label{e:sprF-observables}
\begin{aligned}
\Phi_1(\vec{a}) &= 0.33a_1^2 + 0.27a_1a_2 + 1.28a_1a_3 + 0.88a_2^2+0.49a_2a_3 + 0.05a_3^2, \\
\Phi_2(\vec{a}) &= 0.71a_1^2 + 0.59a_1a_2 + 0.84a_1a_3 + 0.42a_2^2 + 0.83a_2a_3 + 0.31a_3^2, \\
\Phi_3(\vec{a}) &= 0.75a_1^2 + 0.68a_1a_2 + 1.04a_1a_3 + 0.5a_2^2 + 1.52a_2a_3 + 0.38a_3^2, \\
\Phi_4(\vec{a}) &= 0.98a_1^2 + 0.3a_1a_2 + 1.42a_1a_3 + 0.6a_2^2 + 1.21a_2a_3 + 0.02a_3^2.
\end{aligned}
\end{equation}
These were selected from a list of 30 randomly generated quadratic $\Phi$, after removing those whose time average is maximized at one of the equilibria. We deemed this to be the case if the best upper bound on $\overline{\Phi}^*$ obtained with the polynomial optimization techniques of~\cite{Chernyshenko2014,Fantuzzi2016,Goluskin2018,Goluskin2019,Lakshmi2020} differed from the value of $\Phi$ at one of the equilibria by less than 0.01 in absolute value.

\begin{figure}
\centering
\includegraphics[width=1\textwidth]{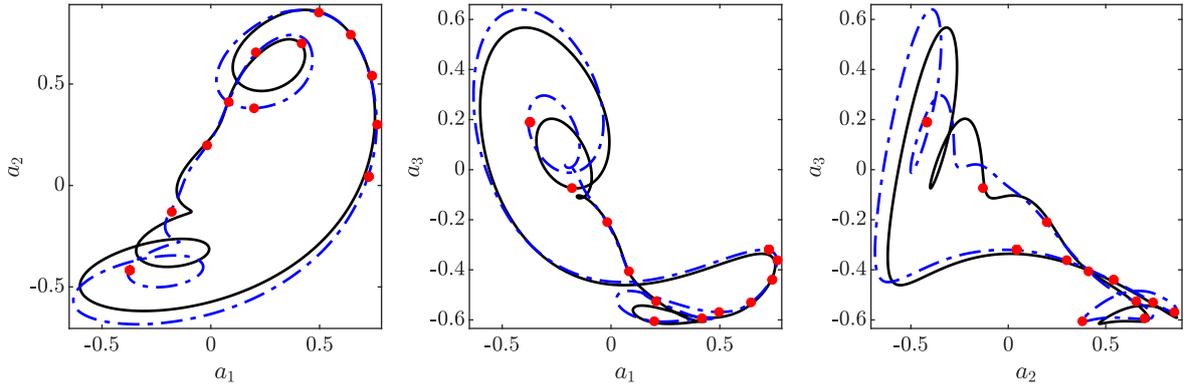}
\caption{Extremal UPO for system~\cref{e:sprF-system} with observable $\Phi_1(\vec{a})$ ({\color{RGBblack}\solidrule}). Also plotted are the converged UPO for the corresponding controlled system~\cref{e:controlled-projection-ode} with $k = 0.25$ ({\color{RGBblue}\dotdashedrule}) and local minimizers of the polynomial $D$ (\mydot{RGBred}) obtained with a degree-14 auxiliary function, which lie inside the set $\mathcal{S}_{10^{-4}}$.}
\label{fig:sprF-Run12-converged-UPOs}
\end{figure}

\begin{figure}
\centering
\includegraphics[width=1\textwidth]{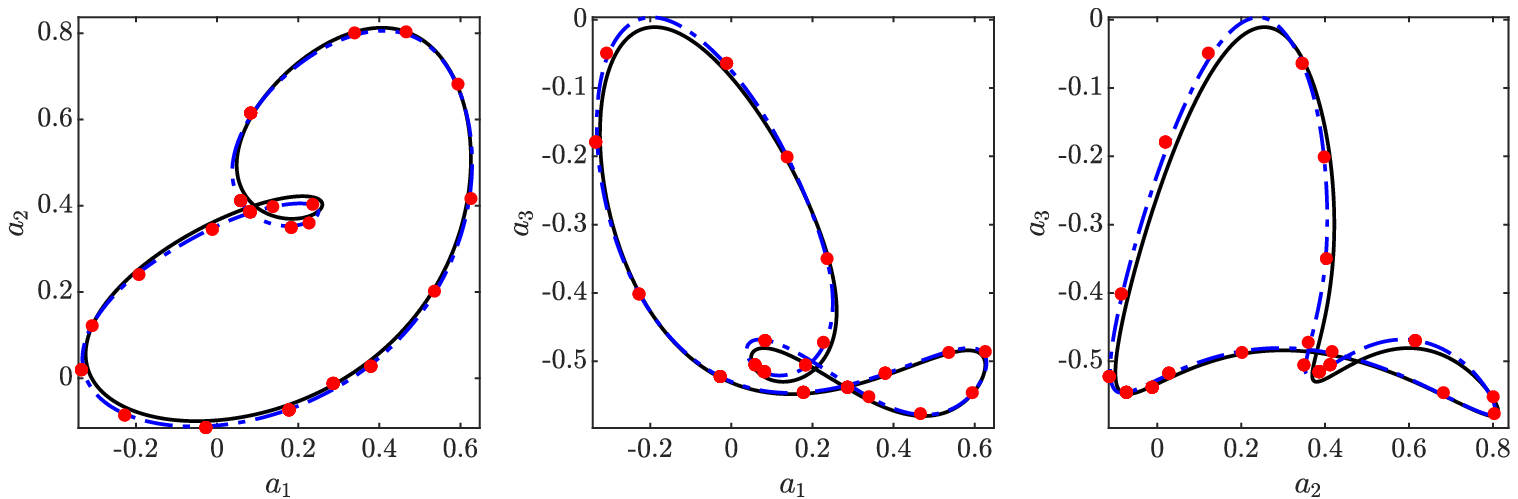}
\caption{Extremal UPO for system~\cref{e:sprF-system} with observable $\Phi_2(\vec{a})$ ({\color{RGBblack}\solidrule}). Also plotted are the converged UPO for the corresponding controlled system~\cref{e:controlled-projection-ode} with $k = 0.25$ ({\color{RGBblue}\dotdashedrule}) and local minimizers of the polynomial $D$ (\mydot{RGBred}) obtained with a degree-14 auxiliary function, which lie inside the set $\mathcal{S}_{10^{-4}}$.}
\label{fig:sprF-Run8-converged-UPOs}
\end{figure}

Results of computations for the observables $\Phi_1$ and $\Phi_2$ are shown in~\Cref{fig:sprF-Run12-converged-UPOs,fig:sprF-Run8-converged-UPOs}. First, we used polynomial optimization to construct a polynomial auxiliary function $V$ of degree 14, which gives the upper bound $\overline{\Phi_1}^* \leq 0.4217$. Increasing the polynomial degree gives the same bound to within 0.047\%, suggesting that our $V$ is very close to being optimal. We then computed local minimizers of the polynomials $D$ corresponding to this near-optimal $V$.  These local minimizers, plotted as red dots in the figure, lie in the set $\mathcal{S}_{10^{-4}}$ and are expected to lie close to the extremal UPO. Recurrence analysis on the trajectory of the controlled ODE system~\cref{e:controlled-projection-ode} with $k=0.25$ starting from the best available local minimizer produced a very good initial guess for a UPO, plotted as dot-dashed blue lines in~\Cref{fig:sprF-Run12-converged-UPOs}, which we converged using the variational algorithm of~\cite{Boghosian2011a}. Numerical continuation of this UPO in $k$ down to $k=0$ resulted in a UPO for~\cref{e:sprF-system}, plotted as a solid black line. The average of $\Phi_1$ over this UPO is within 0.1\% of the upper bound on $\overline{\Phi_1}^*$ reported above, strongly suggesting that we have indeed computed the extremal orbit. Results for the observable $\Phi_2$, plotted in~\Cref{fig:sprF-Run8-converged-UPOs}, lead to a different UPO but are qualitatively analogous, so we do not discuss them for brevity. For both observables, the same results were obtained also when the initial value of the control amplitude $k$ was varied in the range $[0.2,0.7]$. This suggests that, at least for this particular ODE system, our approach is not very sensitive to the initial choice of $k$.


\begin{figure}
\centering
\includegraphics[width=1\textwidth]{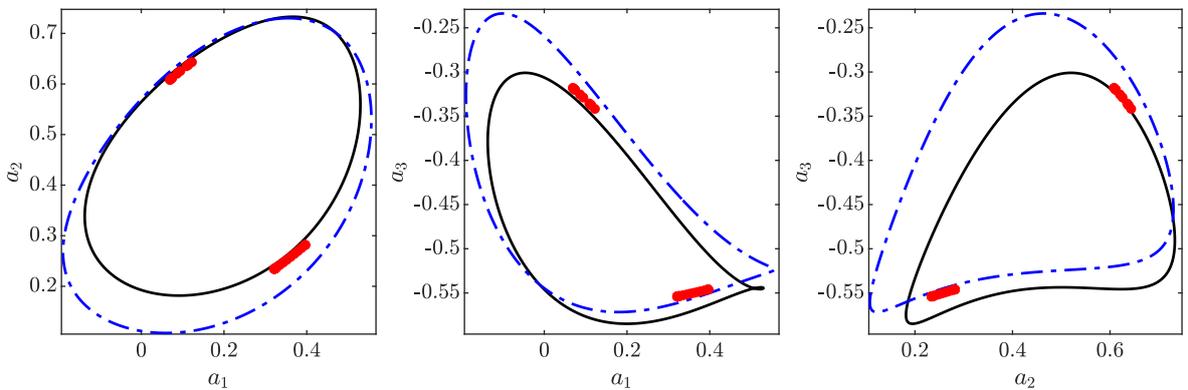}
\caption{
UPO for system~\cref{e:sprF-system} that simultaneously maximizes the time average of $\Phi_3(\vec{a})$ and $\Phi_4(\vec{a})$ ({\color{RGBblack}\solidrule}). Also plotted are the converged UPO for the corresponding controlled system~\cref{e:controlled-projection-ode} with $k = 0.1$ ({\color{RGBblue}\dotdashedrule}) and local minimizers of the polynomial $D$ (\mydot{RGBred}) obtained with a degree-6 auxiliary function, which lie inside the set $\mathcal{S}_{10^{-8}}$.}
\label{fig:sprF-Run1-converged-UPOs}
\end{figure}

Similar calculations for the observables $\Phi_3$ and $\Phi_4$ in~\cref{e:sprF-observables} led to the discovery of only a third UPO, illustrated in~\Cref{fig:sprF-Run1-converged-UPOs}, that simultaneously maximizes the average of both quantities. We conclude this because numerical upper bounds on $\overline{\Phi_3}^*$ and  $\overline{\Phi_4}^*$ computed with degree-10 polynomial auxiliary functions are actually 0.0036\% \textit{smaller} than the averages of $\Phi_3$ and $\Phi_4$ on this UPO.\footnote{This apparent contradiction is due to unavoidably finite tolerances in the algorithms used to optimize the upper bounds, which may return slightly infeasible answers.} However, this orbit could be found using polynomial auxiliary functions of degree as low as six, which yield upper bounds 0.0038\% and 0.0041\% larger than the respective true averages, and with control amplitude $k$ as small as $0.1$. Auxiliary functions of low polynomial degrees are significantly cheaper to construct~\cite{Chernyshenko2014,Fantuzzi2016,Goluskin2018,Goluskin2019,Lakshmi2020} and may be expected to work when the extremal UPO being sought has a simple shape in state space. Indeed, the best control term in~\cref{e:controlled-projection-ode} is obtained when the function $D$ corresponds to an optimal auxiliary function $V$. In this ideal case, $D$ must vanish on the extremal UPO (cf. \cref{s:auxiliary-function-method}). Therefore, a polynomial $V$ with enough degrees of freedom to approximately satisfy the same constraint may lead to good control terms in~\cref{e:controlled-projection-ode}, even if the corresponding upper bound on $\overline{\Phi}^*$ are far from sharp. We suspect that this is why the extremal UPO for $\Phi_3$ and $\Phi_4$ could be found with degree-6 $V$ while the extremal UPOs for $\Phi_1$ and $\Phi_2$, which have a more complicated shape, required increasing the polynomial degree to 14.

A final interesting observation is that, in~\Cref{fig:sprF-Run12-converged-UPOs,fig:sprF-Run8-converged-UPOs}, many of the local minimizers of $D(\vec{a})$ approximating the extremal UPO actually to lie closer to the converged UPO for the \textit{controlled} system. A possible explanation for this is that the extremal UPO for~\cref{e:sprF-system} is not a periodic orbit for the uncontrolled one when $V$ is suboptimal, and need not lie entirely in the set $\mathcal{S}_\varepsilon$ if $\varepsilon$ is small. On the other hand, the local minimizers of $D(\vec{a})$ do lie in $\mathcal{S}_\varepsilon$ by construction, and the control term in~\cref{e:controlled-projection-ode} pushes trajectories exactly towards these points. Consequently, local minimizers of $D$ could be better initial conditions for the controlled system, rather than for the uncontrolled one. This, however, may not always be the case (see, for instance,~\Cref{fig:sprF-Run1-converged-UPOs}) and appears to depend on how suboptimal the auxiliary function is, on the chosen value of the control parameter $k$, and on how ``flat" the polynomial $D$ is in the vicinity of the extremal UPO along directions perpendicular to it. 

\subsection{Lorenz--96 system}
\label{ss:L96-system}
\begin{figure}
\centering
\includegraphics[width=1\textwidth]{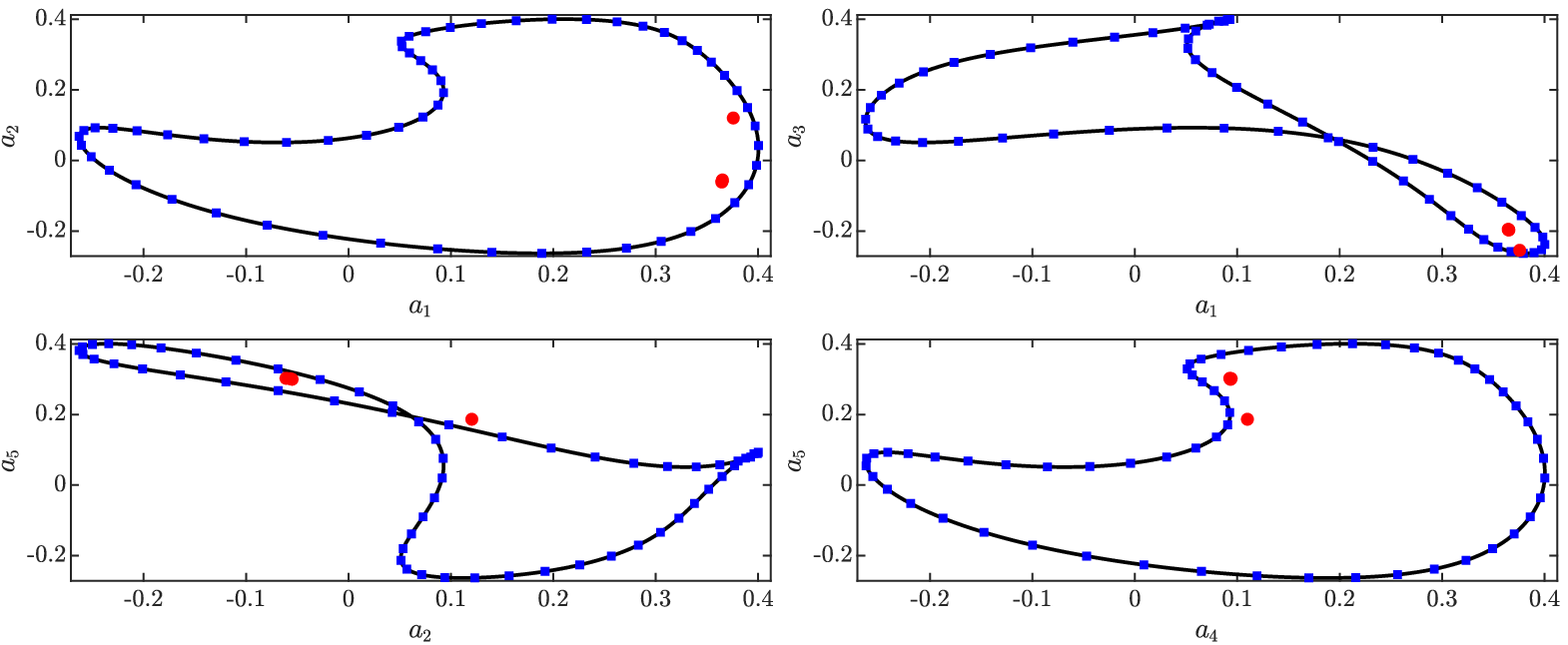}
\caption{UPO for system~\cref{e:L96-system} that maximizes the time average of $\Phi(\vec{a}) = (a_1 - F)^2 + (a_4 - F)^2$ for $F = 8$ ({\color{RGBblack}\solidrule}). Also plotted are the converged UPO for the corresponding controlled system~\cref{e:controlled-projection-ode} with $k = 0.5$ (\mysolidsquare{RGBblue}) and local minimizers of the polynomial $D$ (\mydot{RGBred}) obtained with a degree-10 auxiliary function, which lie inside the set $\mathcal{S}_{10^{-4}}$.}
\label{fig:L96-converged-UPOs}
\end{figure}

We next study the five-dimensional Lorenz--96 ODE system~\cite{Lorenz1996}
\begin{equation}
\label{e:L96-system}
\frac{\mathrm{d} a_i}{\dt} = (a_{i+1} - a_{i-2})a_{i-1} - a_i + F, \qquad i = 1,\dots,5,
\end{equation}
where we adopt the convention that $a_{-1} = a_4$, $a_0 = a_5$ and $a_6 = a_1$. The scalar $F$ is a constant forcing term and the system has a unique equilibrium point $\vec{a}_0 = (F,\dots,F)$. Fixing $F = 8$ we applied the 4-step strategy outlined in~\cref{s:control-suboptimal} with a degree-10 polynomial function and $k = 0.5$ to search for the extremal UPO for the observable
\begin{equation}
    \Phi(\vec{a}) = (a_1 - F)^2 + (a_4 - F)^2,
\end{equation}
whose time average is certainly not maximized by the equilibrium.
Local minima of the polynomial $D$, the converged UPO for the controlled system~\cref{e:controlled-projection-ode}, and the final UPO with $k=0$ are illustrated in~\Cref{fig:L96-converged-UPOs}. The fact that the upper bound on $\overline{\Phi}^*$ obtained in step 1 exceeds the numerical of average of $\Phi$ over the converged UPO by only 0.43\% strongly suggests that we have indeed found the extremal one.

It is instructive to ask whether our control strategy offers any advantages over the method from~\cite{Lakshmi2020}, which attempts to find the extremal UPO directly by using local minimizers of $D$ (supplemented with a guess for the orbit's period) as initial conditions for a single-shooting Newton--Raphson method. With the set of local minimizers shown in~\Cref{fig:L96-converged-UPOs}, the single-shooting Newton--Raphson algorithm fails to converge to any UPO for the Lorenz--96 system~\cref{e:L96-system} even when the initial guess for the period is taken to be equal to the period $T^*$ of the target extremal UPO. To understand this failure in convergence, we integrated both the uncontrolled Lorenz--96 system and the corresponding controlled system with $k = 0.4$ over a time interval of length $T^*$, taking the best local minimizer of $D$ as the initial condition. The shooting error $\|\vec{a}(T^*\ssep \vec{a}_0)- \vec{a}_0\|/\|\vec{a}_0\|$ was 0.256 for the controlled trajectory and 0.396 for the uncontrolled one. This clearly demonstrates the efficacy of the control term in~\cref{e:controlled-projection-ode} in reducing the instability of the orbit, which is too unstable for the shooting strategy of~\cite{Lakshmi2020} to work. This difficulty could not be resolved by working with a more accurate polynomial auxiliary function of degree-14, highlighting that the sensitivity of the single-shooting approach to poor initial conditions cannot be easily mitigated by increasing the accuracy of $V$. In contrast, the control strategy described in this paper works robustly.

\subsection{A model of shear flow}
\label{ss:9-mode-system}
As a final example, we consider a nine-dimensional ODE system modeling sinusoidally forced shear flow in a periodic channel~\cite{Moehlis2004}. The system takes the form
\begin{equation}
\label{e:9-mode-system}
\frac{da_i}{dt} = \lambda_1 \delta_{1i} - \frac{1}{ \Rey} \lambda_{j}a_j + N_{ijk} a_j a_k, \qquad i,j,k = 1,\dots,9,
\end{equation}
where summation over indices $j$ and $k$ is assumed, $\delta_{1i}$ is the usual Kronecker delta, $\Rey$ is a fixed constant representing the flow's Reynolds number, and $\lambda_{j}$, $N_{ijk}$ are numerical coefficients corresponding to the ``NBC'' configuration in~\cite{Moehlis2004,Moehlis2005}.

At all values of $\Rey$, system~\cref{e:9-mode-system} has a locally stable equilibrium point $\vec{a}_l = (1,0,\dots,0)$, which represents the laminar flow state. Trajectories display chaotic behavior as $\Rey$ is raised and a large number of UPOs have been computed for $\Rey>80.54$~\cite{Moehlis2004,Moehlis2005}. Here, we fix $\Rey=120$ and look for the UPO that maximizes (in a time-averaged sense) the energy of perturbations from the laminar flow, 
\begin{equation}
\Phi(\vec{a}) := \|\vec{a} - \vec{a}_l\|^2.
\end{equation}

\begin{figure}[t]
\centering
\includegraphics[width=1\textwidth]{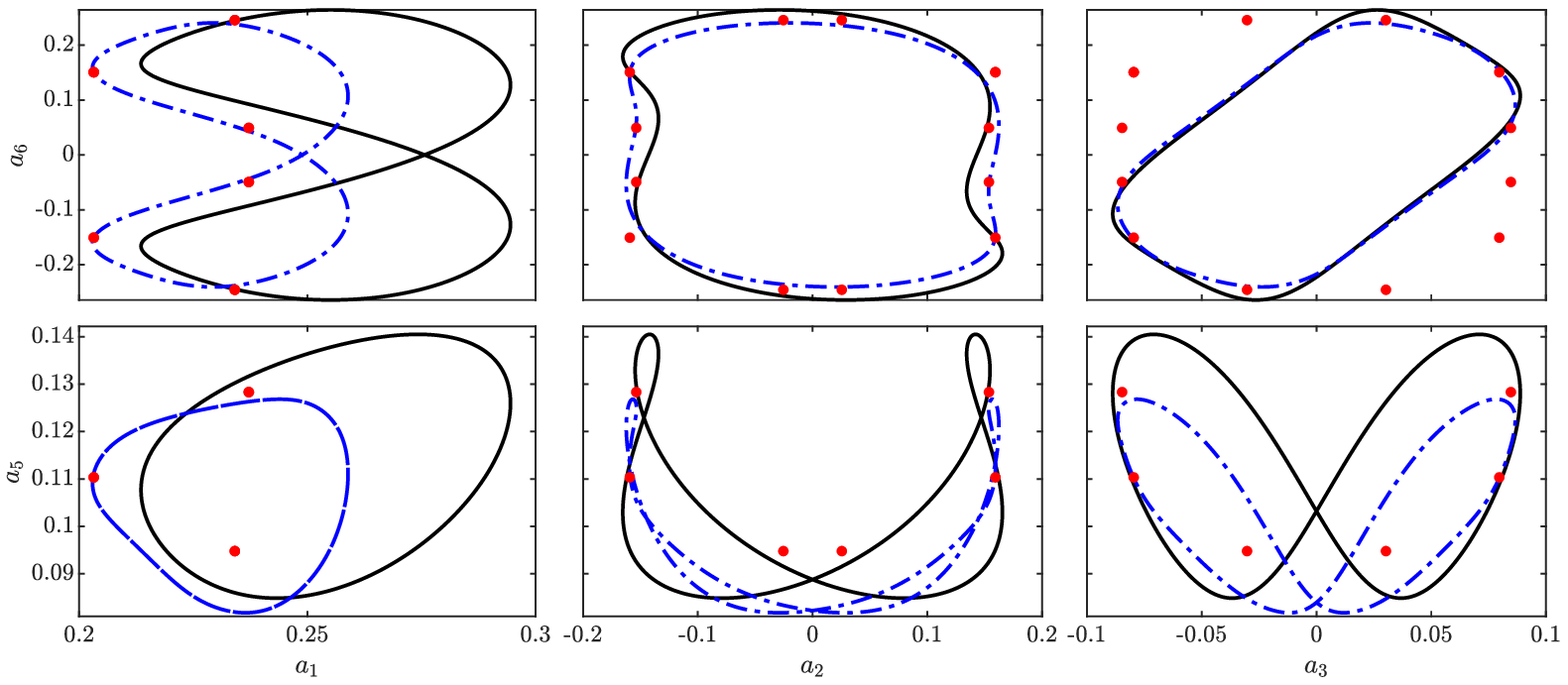}
\caption{
UPO for system~\cref{e:9-mode-system} that maximizes the time average of $\Phi(\vec{a}) = \|\vec{a} - \vec{a}_l\|^2$ at $\Rey = 120$ ({\color{RGBblack}\solidrule}). Also plotted are the converged UPO for the corresponding controlled system~\cref{e:controlled-projection-ode} with $k = 0.4$ ({\color{RGBblue}\dotdashedrule}) and local minimizers of the polynomial $D$ (\mydot{RGBred}) obtained with a degree-10 auxiliary function, which lie inside the set $\mathcal{S}_{10^{-4}}$.
}
\label{fig:9Mode-converged-UPOs}
\end{figure}

The result of computations with a degree-10 auxiliary function and $k = 0.4$, are shown in~\Cref{fig:9Mode-converged-UPOs}. The numerical average of $\Phi$ over the converged UPO for system~\cref{e:9-mode-system} is only 0.63\% less than our upper bound on $\Phi$. Once again, this strongly suggests that our method yields the extremal UPO for $\Phi$.
As in~\cref{ss:L96-system}, the uncontrolled single-shooting method of~\cite{Lakshmi2020} fails to produce any periodic orbit, even when the initial guess for the period is taken to be equal to the period $T^*$ of the target extremal UPO. Integrating both the uncontrolled system~\cref{e:9-mode-system} and its controlled counterpart~\cref{e:controlled-projection-ode} with $k = 0.4$ over a time interval of length $T^*$ using the best local minimizer of $D$ as the initial condition, we find that the shooting error $\|\vec{a}(T^*\ssep \vec{a}_0)- \vec{a}_0\|/\|\vec{a}_0\|$ is 0.0174 for the controlled trajectory and 0.3751 for the uncontrolled one. This again illustrates the stabilizing effect of the control. Of course, we cannot exclude that the single-shooting method of~\cite{Lakshmi2020} will work if one uses polynomial auxiliary functions of higher degree, but doing so would require significantly larger computational resources than those available to this study. The control strategy approach presented here, however, could be implemented successfully without difficulties.

\section{Discussion}
\label{s:discussion}
The examples above demonstrate that the four-step strategy introduced in this work enables one to compute extremal UPOs robustly, at least for low-dimensional ODE systems. The ability to focus on extremal orbits is a particular advantage of our method, because performing recurrence analysis on chaotic trajectories might easily miss UPOs with extreme behavior that do not contribute significantly to the statistics of the chaotic attractor. 

On the other hand, the need to construct polynomial auxiliary functions that produce sufficiently accurate upper bounds on the extremal time average $\overline{\Phi}^*$ of the quantity of interest currently poses a significant barrier to scalability. This is because the computational complexity of the polynomial optimization techniques described in~\cite{Chernyshenko2014,Fantuzzi2016,Goluskin2018,Goluskin2019,Lakshmi2020} grows combinatorially as the number of ODE states and the degree of the auxiliary function is increased~\cite{Lasserre2015}. As a result, the largest polynomial degree that can currently be considered for the nine-dimensional system in \cref{ss:9-mode-system} is approximately 10 on a workstation with 64GB of RAM, and reduces to no more than 4 or 6 for ODEs with a few tens of states. Nevertheless, removing computational bottlenecks in general polynomial optimization and in its applications to dynamical systems are problems that have attracted significant interest in recent years (see, e.g.,~\cite{Zheng2017csl,Lasserre2017b,Weisser2017,Ahmadi2018,Ahmadi2019,Zheng2019sos,Tacchi2019,Wang2020a,Wang2020b,Schlosser2020}), so we expect that our UPO search strategy will become practical for ODEs of moderate dimension in the near future.

With computational aspects in mind, a particularly attractive aspect of the control strategy proposed in this paper is that its four steps do not depend on the particular algorithms used to carry them out. This enables one to use more sophisticated numerical techniques not only to construct polynomial auxiliary functions, but also to converge UPOs for the controlled system~\cref{e:controlled-projection-ode}. For example, one could augment multiple-shooting and variational algorithms with recent approaches to identifying near-periodic trajectory segments based on dynamic mode decomposition~\cite{Page2020}, which are more robust than the traditional recurrence analysis employed in this work.

Finally, all results reported in this paper were obtained by fixing the control amplitude $k$ in~\cref{e:controlled-projection-ode} to an arbitrary constant value in the interval $(0,1)$. Methods to optimize $k$ statically or dynamically by taking into account the characteristics of the polynomial $D(\vec{a})$ would be a highly valuable addition to our approach, because simply fixing a large $k$ value would mean that the periodic orbit for the controlled system~\cref{e:controlled-projection-ode}, if it exists, lies far away from the extremal UPO for the uncontrolled system~\cref{e:ode}. This situation is not ideal, because then a continuous branch connecting a periodic orbit for the controlled system to the extremal UPO for the uncontrolled one may not exist. On the other hand, taking $k$ too small could mean that the stabilizing effect of the control term is not sufficient to find any near-periodic trajectory segments. There clearly is an optimal choice for $k$, but there are also many ways to formulate an optimal control problem for the controlled ODEs~\cref{e:controlled-ode} or~\cref{e:controlled-projection-ode} and, without \textit{a priori} knowledge of an extremal UPO and/or of the polynomial $D$, it is not immediately clear which one leads to the best UPO approximation. We leave further investigation of this problem to future work.

\section{Conclusion} 
\label{s:conclusion}
We have presented a novel technique of computing UPOs for ODE systems governed by polynomial dynamics, which overcomes the deficiencies of a related approach presented in~\cite{Lakshmi2020}. As in that work, the UPO search procedure is initiated by leveraging polynomial optimization techniques to construct an auxiliary function $V$ that proves a near-sharp bound $U$ on the maximal value of the infinite-time average of a prescribed observable $\Phi$. Direct unconstrained minimization of the polynomial $D = U - \Phi - \vec{f} \cdot \nabla V$ then yields a set of points which lie close to the extremal UPO for the original system. The novel contribution of this work was to show that this polynomial $D$ can also be used to construct an effective control strategy which reduces the instability of the orbit, and aids its computation with traditional techniques. 

More precisely, we have formulated a family of controlled ODE systems~\cref{e:controlled-projection-ode} parameterized by a control amplitude $k$, which reduce to the original ODE system when $k = 0$. When the auxiliary function $V$ is optimal and the corresponding $D$ satisfies the assumptions of~\Cref{p:proposition}, a sufficiently large control amplitude $k$ will guarantee the existence of a (locally) stable periodic orbit for the controlled system that can be found simply by time integration. By construction, this periodic orbit coincides with the UPO of the original, uncontrolled system that maximizes the infinite-time average of the given observable $\Phi$. If $V$ is not optimal, which is often the case in practice, stabilization cannot be guaranteed. However, for near-optimal $V$ one expects that periodic orbits for the controlled systems not only still exist, but are also less unstable than the extremal UPOs for the original ODE and, crucially, connect to it continuously as the control amplitude $k$ is reduced to zero. One can therefore first compute an orbit for the controlled system by combining recurrence analysis along trajectories starting from local minimizers of $D$ with shooting or variational algorithms, and then continue it numerically in $k$ to recover the extremal UPO for the original ODE, which is harder to compute directly.

This process was applied successfully in~\cref{s:results} to a three-dimensional system with a chaotic attractor \cite{Sprott1994}, a five-dimensional version of the Lorenz--96 system~\cite{Lorenz1996}, and a nine-dimensional model of shear flow~\cite{Moehlis2004}. The control methodology developed in this paper was essential to compute the correct extremal UPOs for a number of different observables $\Phi$. In stark contrast, the uncontrolled approach in~\cite{Lakshmi2020} failed in all cases. We cannot say if the robust numerical behavior observed in these examples is generic, and further theoretical and computational advances are necessary before complex high-dimensional ODE systems can be tackled using the ideas we have described here. Nevertheless, augmenting traditional techniques for converging UPOs with recent frameworks for ODE analysis via polynomial optimization promises to be a fruitful avenue of research.

\pdfbookmark{Acknowledgments}{bookmark:Funding}
\paragraph*{Acknowledgments.}
MVL was supported by EPSRC studentship award 2092930 under grant EP/N509486/1. GF was supported by an Imperial College Research Fellowship.

\pdfbookmark{Appendices}{bookmark:Appendices}
\appendix

\section{A stability result for ODE \texorpdfstring{\Cref{e:controlled-projection-ode}}{(\ref{e:controlled-projection-ode})}}
\label{app:stability-result}

In the case of optimal $V$, the stability argument in~\Cref{p:proposition} can be extended to the controlled system~\cref{e:controlled-projection-ode} under moderate assumptions on $D$ and $\vec{f}$. Proceeding as in the proof of~\Cref{p:proposition}, one must show that the right-hand side of the inequality
\begin{equation}
\frac{\mathrm{d}D}{\dt}
\leq \max_{\vec{a}\in \mathcal{V}_\gamma}\abs{\vec{f} \cdot \nabla D} - k \|\nabla D \|^2 \left(1- \frac{\vert\vec{f} \cdot \nabla D\vert^2}{\|\vec{f}\|^2 \, \|\nabla D \|^2} \right)
\end{equation}
is negative at all points $\vec{a} \in \mathcal{V}_\gamma$ for sufficiently large $k$. This is true if $D$ and $\vec{f}$ satisfy
\begin{equation}\label{e:necessary-condition-controlled-stability}
    \frac{\vert\vec{f} \cdot \nabla D\vert}{\|\vec{f}\| \, \|\nabla D \|} < \text{const} < 1 \qquad \forall \vec{a} \in \mathcal{V}_\gamma.
\end{equation}
This constraint is weak, because
$\nabla (\vec f\cdot \nabla D) = \nabla \vec{f} \ \nabla D + \vec{f}\cdot (\nabla \otimes \nabla D)$ vanishes on the orbit: the first term in the sum is zero there because $\nabla D = \vec{0}$, while the second term is zero because it is the material derivative of $\nabla D$ along the orbit. As a result, one can expect that in the vicinity of the orbit $\vec{f} \cdot \nabla D\sim d^2$ while $\|\nabla D \| \sim d$, where $d$ is the distance to the orbit. This, in particular, will be the case at all the points on the orbit where the Hessian matrix $\nabla \otimes \nabla D$ is positive definite in the subspace perpendicular to $\vec{f}$. One can therefore ensure~\cref{e:necessary-condition-controlled-stability} by taking $\gamma$ small enough.

\section{Numerical implementation}
\label{ss:numerical-implementation}
For all examples presented in \cref{s:results}, polynomial auxiliary functions were constructed with the polynomial optimization framework described in~\cite{Chernyshenko2014,Fantuzzi2016,Goluskin2018,Goluskin2019,Lakshmi2020} using a customized version\footnote{Available from \url{https://github.com/aeroimperial-optimisation/aeroimperial-yalmip}} of the MATLAB optimization toolbox YALMIP~\cite{Lofberg2004} and optimization solver MOSEK v.9.2~\cite{MOSEKApS2020}. Approximate local minima of the corresponding polynomial $D(\vec{a})$ were found using the BFGS quasi-Newton method~\cite{Broyden1970,Fletcher1970,Goldfarb1970,Shanno1970} implemented in MATLAB's built-in function \texttt{fminunc}, with the step tolerance (relative lower bound on the size of an iteration step) and the first-order optimality tolerance (lower bound on $\| \nabla D\|_{\infty}$) both set to $10^{-16}$.

The best local minimum was then used to integrate the controlled system~\cref{e:controlled-projection-ode} forward in time for a value of $k$ fixed arbitrarily, and recurrence analysis was employed to identify near-periodic trajectory segments~\cite{Viswanath2007,Cvitanovic2010}. In our implementation, a near recurrence of period $T$ was deemed to have occurred at time $t$ if the quantity
\begin{equation}
\label{e:recurrence-function}
R(t,T) := \frac{ \| \vec{a}(t) - \vec{a}(t-T)\| }{ \| \vec{a}(t) \| } \leq 0.025.
\end{equation}
The portion of trajectory between times $t-T$ and $T$ was then used as an initial guess to converge a periodic orbit using the variational approach described in~\cite{Boghosian2011a}. This method requires minimizing the cost function
\begin{equation}
\label{e:cost-function}
C(\hat{\vec{a}}_0, \dots, \hat{\vec{a}}_{N-1}, T) := 
\frac{N}{2T}
\sum_{i=0}^{N-1} 
\left\|  
\hat{\vec{a}}_{i+1} - \hat{\vec{a}}_i -
{\frac{T}{N}} \, \vec{h}_k\!\left(\frac{\hat{\vec{a}}_{i+1} + \hat{\vec{a}}_i}{2}\right)
\right\|^2
\end{equation}
over the orbit's period $T$ and over $N$ points $\hat{\vec{a}}_0, \hat{\vec{a}}_1, \ldots, \hat{\vec{a}}_{N-1}$ distributed at equal time intervals $T/N$ along the orbit. (In writing~\cref{e:cost-function}, we used the convention that $\hat{\vec{a}}_{N}=\hat{\vec{a}}_0$ by periodicity.) We solved this nonlinear least-squares problem using the Levenberg--Marquardt algorithm~\cite{Levenberg1944,Marquardt1963} implemented in the MATLAB built-in function \texttt{lsqnonlin}, with the function tolerance on $C$ (relative lower bound on the change in the value of $C$ after an iteration) and the step tolerance on the vector $(\hat{\vec{a}}_0, \dots, \hat{\vec{a}}_{N-1},T)$ both set to $10^{-16}$.

Periodic orbits for the controlled ODE system~\cref{e:controlled-projection-ode} computed in this way were then numerically continued in $k$ down to $k=0$. This could be done with a variety of sophisticated numerical approaches (see, e.g.,~\cite{Dhooge2008}); here, however, we simply minimized~\cref{e:cost-function} at increasingly small values of $k$ using the minimizer from the previous computation as the initial condition.

\pdfbookmark{References}{bookmark:References}
\bibliography{./21-01-22-arXiv-submission.bib}

\begin{thebibliography}{10}

\bibitem{Cvitanovic1991}
P.~Cvitanovi{\'c}.
\newblock Periodic orbits as the skeleton of classical and quantum chaos.
\newblock {\em Phys. D}, 51:138--151, 1991.

\bibitem{Artuso1990}
R.~Artuso, E.~Aurell, and P.~Cvitanovic.
\newblock Recycling of strange sets: {{I}}. {{Cycle}} expansions.
\newblock {\em Nonlinearity}, 3:325--359, 1990.

\bibitem{Cvitanovic1995}
P.~Cvitanovi{\'c}.
\newblock Dynamical averaging in terms of periodic orbits.
\newblock {\em Phys. D}, 83:109--123, 1995.

\bibitem{Kelley2003}
C.~T. Kelley.
\newblock {\em Solving {{Nonlinear Equations}} with {{Newton}}'s {{Method}}}.
\newblock Fundam. {{Algorithms}}. {SIAM}, {Philadelphia}, 2003.

\bibitem{Ascher1995a}
U.~M. Ascher, R.~M.~M. Mattheij, and R.~D. Russell.
\newblock Initial {{Value Methods}}.
\newblock In {\em Numerical {{Solution}} of {{Boundary Value Problems}} for
  {{Ordinary Differential Equations}}}, Classics {{Appl}}. {{Math}}.,
  chapter~4, pages 132--184. {SIAM}, {Philadelphia}, 1995.

\bibitem{Cvitanovic2017}
P.~Cvitanovi{\'c}, R.~Artuso, R.~Mainieri, G.~Tanner, and G.~Vattay.
\newblock {\em Chaos: {{Classical}} and {{Quantum}}}.
\newblock {Niels Bohr Institute}, {Copenhagen}, 2017.
\newblock (version 15.9).

\bibitem{Lan2004}
Y.~Lan and P.~Cvitanovi{\'c}.
\newblock Variational method for finding periodic orbits in a general flow.
\newblock {\em Phys. Rev. E}, 69:016217, 2004.

\bibitem{Boghosian2011a}
B.~M. Boghosian, L.~M. Fazendeiro, J.~L{\"a}tt, H.~Tang, and P.~V. Coveney.
\newblock New variational principles for locating periodic orbits of
  differential equations.
\newblock {\em Philos. Trans. Roy. Soc. A}, 369:2211--2218, 2011.

\bibitem{Fazendeiro2010}
L.~Fazendeiro, B.~M. Boghosian, P.~V. Coveney, and J.~L{\"a}tt.
\newblock Unstable periodic orbits in weak turbulence.
\newblock {\em J. Comput. Sci.}, 1:13--23, 2010.

\bibitem{Chandler2013}
G.~J. Chandler and R.~R. Kerswell.
\newblock {Invariant recurrent solutions embedded in a turbulent
  two-dimensional Kolmogorov flow}.
\newblock {\em J. Fluid Mech.}, 722:554--595, 2013.

\bibitem{Lucas2014}
D.~Lucas and R.~Kerswell.
\newblock {Spatiotemporal dynamics in two-dimensional Kolmogorov flow over
  large domains}.
\newblock {\em J. Fluid Mech.}, 750:518--554, 2014.

\bibitem{Chernyshenko2014}
S.~I. Chernyshenko, P.~Goulart, D.~Huang, and A.~Papachristodoulou.
\newblock Polynomial sum of squares in fluid dynamics: A review with a look
  ahead.
\newblock {\em Philos. Trans. Roy. Soc. A}, 372:20130350, 2014.

\bibitem{Fantuzzi2016}
G.~Fantuzzi, D.~Goluskin, D.~Huang, and S.~I. Chernyshenko.
\newblock Bounds for {{Deterministic}} and {{Stochastic Dynamical Systems}}
  using {{Sum}}-of-{{Squares Optimization}}.
\newblock {\em SIAM J. Appl. Dyn. Syst.}, 15:1962--1988, 2016.

\bibitem{Goluskin2018}
D.~Goluskin.
\newblock Bounding {{Averages Rigorously Using Semidefinite~Programming}}:
  {{Mean Moments}} of the {{Lorenz System}}.
\newblock {\em J. Nonlinear Sci.}, 28:621--651, 2018.

\bibitem{Tobasco2018}
I.~Tobasco, D.~Goluskin, and C.~R. Doering.
\newblock Optimal bounds and extremal trajectories for time averages in
  nonlinear dynamical systems.
\newblock {\em Phys. Lett. A}, 382:382--386, 2018.

\bibitem{Korda2021}
M.~Korda, D.~Henrion, and I.~Mezi{\'{c}}.
\newblock {Convex computation of extremal invariant measures of nonlinear
  dynamical systems and Markov processes}.
\newblock {\em J. Nonlinear Sci.}, 31:14(1--26), 2021.

\bibitem{Yang2000}
T.-H. Yang, B.~R. Hunt, and E.~Ott.
\newblock Optimal periodic orbits of continuous time chaotic systems.
\newblock {\em Phys. Rev. E}, 62:1950--1959, 2000.

\bibitem{Lakshmi2020}
M.~V. Lakshmi, G.~Fantuzzi, J.~D. {Fern{\'a}ndez-Caballero}, Y.~Hwang, and
  S.~I. Chernyshenko.
\newblock Finding extremal periodic orbits with polynomial optimization, with
  application to a nine-mode model of shear flow.
\newblock {\em SIAM J. Appl. Dyn. Syst.}, 19:763--787, 2020.

\bibitem{Goluskin2019}
D.~Goluskin and G.~Fantuzzi.
\newblock Bounds on mean energy in the {{Kuramoto}}\textendash{{Sivashinsky}}
  equation computed using semidefinite programming.
\newblock {\em Nonlinearity}, 32:1705--1730, 2019.

\bibitem{Lasagna2016}
D.~Lasagna, D.~Huang, O.~R. Tutty, and S.~Chernyshenko.
\newblock Sum-of-squares approach to feedback control of laminar wake flows.
\newblock {\em J. Fluid Mech.}, 809:628--663, 2016.

\bibitem{Parrilo2003}
P.~A. Parrilo.
\newblock Semidefinite programming relaxations for semialgebraic problems.
\newblock {\em Math. Program., Ser. B}, 96:293--320, 2003.

\bibitem{Parrilo2012}
P.~Parrilo.
\newblock Polynomial {{Optimization}}, {{Sums}} of {{Squares}}, and
  {{Applications}}.
\newblock In {\em Semidefinite {{Optimization}} and {{Convex Algebraic
  Geometry}}}, {{MOS}}-{{SIAM Ser}}. {{Optim}}., chapter~3, pages 47--157.
  {SIAM}, {Philadelphia}, 2012.

\bibitem{Lasserre2015}
J.~B. Lasserre.
\newblock {\em An {{Introduction}} to {{Polynomial}} and {{Semi}}-{{Algebraic
  Optimization}}}.
\newblock Cambridge {{Texts Appl}}. {{Math}}. {Cambridge University Press},
  {Cambridge}, 2015.

\bibitem{Sprott1994}
J.~C. Sprott.
\newblock Some simple chaotic flows.
\newblock {\em Phys. Rev. E}, 50:R647--R650, 1994.

\bibitem{Lorenz1996}
E.~N. Lorenz.
\newblock Predictability: A problem partly solved.
\newblock In {\em Proceedings of the {{Seminar}} on {{Predictability}}},
  volume~1 of {\em {{ECWF Seminar}}}, pages 1--18, {Reading, UK}, 1996.
  {ECMWF}.

\bibitem{Moehlis2004}
J.~Moehlis, H.~Faisst, and B.~Eckhardt.
\newblock A low-dimensional model for turbulent shear flows.
\newblock {\em New J. Phys.}, 6:56, 2004.

\bibitem{Moehlis2005}
J.~Moehlis, H.~Faisst, and B.~Eckhardt.
\newblock Periodic {{Orbits}} and {{Chaotic Sets}} in a {{Low}}-{{Dimensional
  Model}} for {{Shear Flows}}.
\newblock {\em SIAM J. Appl. Dyn. Syst.}, 4:352--376, 2005.

\bibitem{Zheng2017csl}
Y.~Zheng, G.~Fantuzzi, and A.~Papachristodoulou.
\newblock {Exploiting sparsity in the coefficient matching conditions in
  sum-of-squares programming using ADMM}.
\newblock {\em IEEE Control Syst. Lett.}, 1(1):80--85, 2017.

\bibitem{Lasserre2017b}
J.~B. Lasserre, K.-C. Toh, and S.~Yang.
\newblock {A bounded degree SOS hierarchy for polynomial optimization}.
\newblock {\em EURO J. Comput. Optim.}, 5(1-2):87--117, 2017.

\bibitem{Weisser2017}
T.~Weisser, J.-B. Lasserre, and K.-C. Toh.
\newblock {Sparse-BSOS: a bounded degree SOS hierarchy for large scale
  polynomial optimization with sparsity}.
\newblock {\em Math. Program. Comput.}, pages 1--32, 2017.

\bibitem{Ahmadi2018}
A.~A. Ahmadi, G.~Hall, A.~Papachristodoulou, J.~Saunderson, and Y.~Zheng.
\newblock {Improving efficiency and scalability of sum of squares optimization:
  Recent advances and limitations}.
\newblock In {\em Proceedings of the 56th Annual Conference on Decision and
  Control}, pages 453--462, 2018.

\bibitem{Ahmadi2019}
A.~Ahmadi and A.~Majumdar.
\newblock {{DSOS}} and {{SDSOS Optimization}}: {{More Tractable Alternatives}}
  to {{Sum}} of {{Squares}} and {{Semidefinite Optimization}}.
\newblock {\em SIAM J. Appl. Algebra Geom.}, 3:193--230, 2019.

\bibitem{Zheng2019sos}
Y.~Zheng, G.~Fantuzzi, and A.~Papachristodoulou.
\newblock {Fast ADMM for sum-of-squares programs using partial orthogonality}.
\newblock {\em IEEE Trans. Automat. Control}, 64(9):3869--3876, 2019.

\bibitem{Tacchi2019}
M.~Tacchi, C.~Cardozo, D.~Henrion, and J.-B. Lasserre.
\newblock {Approximating regions of attraction of a sparse polynomial
  differential system}.
\newblock \href{http://arxiv.org/abs/1911.09500}{arXiv:1911.09500} [math.OC],
  2019.

\bibitem{Wang2020a}
J.~Wang, V.~Magron, and J.-B. Lasserre.
\newblock {TSSOS: A moment-SOS hierarchy that exploits term sparsity}.
\newblock {\em SIAM J. Optim.}, 31(1):30--58, 2021.

\bibitem{Wang2020b}
J.~Wang, V.~Magron, and J.-B. Lasserre.
\newblock {Chordal-TSSOS: A moment-SOS hierarchy that exploits term sparsity
  with chordal extension}.
\newblock {\em SIAM J. Optim.}, 31(1):114--141, 2021.

\bibitem{Schlosser2020}
C.~Schlosser and M.~Korda.
\newblock {Sparse moment-sum-of-squares relaxations for nonlinear dynamical
  systems with guaranteed convergence}.
\newblock \href{http://arxiv.org/abs/2012.05572}{arXiv:2012.05572} [math.OC],
  2020.

\bibitem{Page2020}
J.~Page and R.~R. Kerswell.
\newblock Searching turbulence for periodic orbits with dynamic mode
  decomposition.
\newblock {\em J. Fluid Mech.}, 886, 2020.

\bibitem{Lofberg2004}
J.~L{\"o}fberg.
\newblock {{YALMIP}} : A toolbox for modeling and optimization in {{MATLAB}}.
\newblock In {\em 2004 {{IEEE International Conference}} on {{Robotics}} and
  {{Automation}} ({{IEEE Cat}}. {{No}}.{{04CH37508}})}, pages 284--289, {New
  Orleans}, 2004.

\bibitem{MOSEKApS2020}
{MOSEK ApS}.
\newblock {{MOSEK Optimization Toolbox}} for {{MATLAB}}, 2020.
\newblock (version 9.2.28).

\bibitem{Broyden1970}
C.~G. Broyden.
\newblock The {{Convergence}} of a {{Class}} of {{Double}}-rank {{Minimization
  Algorithms}} 1. {{General Considerations}}.
\newblock {\em IMA J. Appl. Math.}, 6:76--90, 1970.

\bibitem{Fletcher1970}
R.~Fletcher.
\newblock A new approach to variable metric algorithms.
\newblock {\em Comput. J.}, 13:317--322, 1970.

\bibitem{Goldfarb1970}
D.~Goldfarb.
\newblock A {{Family}} of {{Variable}}-{{Metric Methods Derived}} by
  {{Variational Means}}.
\newblock {\em Math. Comp.}, 24:23--26, 1970.

\bibitem{Shanno1970}
D.~F. Shanno.
\newblock Conditioning of {{Quasi}}-{{Newton Methods}} for {{Function
  Minimization}}.
\newblock {\em Math. Comp.}, 24:647--656, 1970.

\bibitem{Viswanath2007}
D.~Viswanath.
\newblock Recurrent motions within plane {{Couette}} turbulence.
\newblock {\em J. Fluid Mech.}, 580:339--358, 2007.

\bibitem{Cvitanovic2010}
P.~Cvitanovi{\'c} and J.~F. Gibson.
\newblock Geometry of the turbulence in wall-bounded shear flows: Periodic
  orbits.
\newblock {\em Phys. Scr.}, T142:014007, 2010.

\bibitem{Levenberg1944}
K.~Levenberg.
\newblock A method for the solution of certain non-linear problems in least
  squares.
\newblock {\em Quart. Appl. Math.}, 2:164--168, 1944.

\bibitem{Marquardt1963}
D.~W. Marquardt.
\newblock An {{Algorithm}} for {{Least}}-{{Squares Estimation}} of {{Nonlinear
  Parameters}}.
\newblock {\em SIAM J. Appl. Math.}, 11:431--441, 1963.

\bibitem{Dhooge2008}
A.~Dhooge, W.~Govaerts, Y.~A. Kuznetsov, H.~G. Meijer, and B.~Sautois.
\newblock New features of the software {{MatCont}} for bifurcation analysis of
  dynamical systems.
\newblock {\em Math. Comput. Model. Dyn. Syst.}, 14:147--175, 2008.

\end{thebibliography}
\end{document}